 
\mag=\magstep1
\documentstyle{amsppt}
\input amsppt1
\pageheight{22true cm}
\pagewidth{16.5true cm}
\parindent=4mm
\parskip=3pt plus1pt minus.5pt
\nologo\NoRunningHeads\NoBlackBoxes

\def\i{\looparrowright}
\def\e{\hookrightarrow}

\def\f{\flushpar }
\def\nl{\newline }
\def\np{\newpage }
\def\x{\times }

\document

\topmatter
\title
Singularities of the projections of $n$-dimensional knots
\endtitle
\author
Eiji Ogasa
\endauthor
\thanks{
This research was partially supported by Research Fellowships 
of the Promotion of Science for Young Scientists.}
\affil
Department of Mathematical Sciences\\
 University of Tokyo, Komaba\\
 Tokyo 153,   Japan\\    
pqr100pqr100\@yahoo.co.jp   \\
This paper is published in\\ 
Mathematical Proceedings of the Cambridge Philosophical Society,  126, 1999, 511-519.\\ 
This manuscript is not the published version.\\
\endaffil
\endtopmatter

\baselineskip11pt
\f{\bf Abstract.}   
Let $n\geqq5$. 
There is a smoothly knotted $n$-dimensional sphere in $(n+2)$-space 
such that the singular point set of its projection in $(n+1)$-space 
 consists of double points
and that 
the components of the singular point set are two.
(The sphere is {\it knotted } in the sense that
 it does not bound any embedded $(n+1)$-ball in $(n+2)$-space.) 
Furthermore, the projection is not the projection of any unknotted sphere in $(n+2)$-space. 
There are two inequivalent embeddings of an $n$-manifold in $(n+2)$-space such that 
the projection of one of these in $(n+1)$-space has no double points and 
the projection of the other has a connected embedded double point set.

\np
\head 1.Introduction and Main results \endhead     

In the study of classical knots,    
the projections of a knot into  $\Bbb R^2$ plays an important role
(see e.g. 
[A], 
[BL], 
[C], 
[CF], 
[J], 
[Kf1], 
[Ko], 
[Re], 
[V], 
[W], 
for example).   
For 2-dimensional knots in $\Bbb R^4$,   
 the projection in  $\Bbb R^3$ is considered 
(see 
[CS2], 
[Km], 
[KSS], 
[Su], 
for example.).  
Projections of codimension two submanifolds of $\Bbb R^n$ into $\Bbb R^{n-1}$ 
is studied by 
[CS1], 
[CS3], 
[G], 
[R], 
etc.
Projections of $p$-dimensional submanifolds of $\Bbb R^n$ into $\Bbb R^{n-1}$,  
for $p<n-2$, 
is studied by 
[SS], 
etc.

In this paper, we consider the projection (into $\Bbb R^{n+1}$) of $n$-dimensional embeddings in $\Bbb R^{n+2}$ where $n\geqq5$. 
We work in the smooth category throughout. 
We consider those embeddings for which the projection has relatively simple self-intersections. 
We show that there are embeddings that are truly knotted, but whose projections have simple 
self intersections. First we introduce some notation.

We work in the smooth category.

An {\it (oriented)  n-(dimensional) knot $K$} is 
a smooth oriented submanifold of $\Bbb R^{n+1}$ $\x$ $\Bbb R$ 
which is diffeomorphic to the standard $n$-sphere.   
We say that n-knots $K_1$ and $K_2$ are {\it equivalent} 
if there exists an orientation preserving diffeomorphism 
$f:$ $\Bbb R^{n+1}$ $\x$ $\Bbb R$ $\rightarrow$ $\Bbb R^{n+1}$ $\x$ $\Bbb R$
such that $f(K_1)$=$K_2$  and 
$f\vert_{K_1}:$ $K_1$ $\rightarrow$ $K_2$ is 
an orientation preserving diffeomorphism.    

Note. In many other papers, including the author's, the definition of $n$-knot is  
a smooth oriented submanifold of $\Bbb R^{n+1}$ $\x$ $\Bbb R$ 
which is PL homeomorphic to the standard $n$-sphere.   
But in this paper, we adopt the former one and reject the latter one.

Let $T$ be the unit $n$-sphere of 
$\Bbb R^{n+1}$ $\x$ $\{0\}$ $\subset$ $\Bbb R^{n+1}$ $\x$ $\Bbb R$.   
Then $T$ is an $n$-knot. 
An $n$-knot $K$ is said to be {\it unknotted} if $K$ is equivalent to $T$.  
If $K$ is not unknotted, then we say that $K$ is {\it truly knotted}.

Let 
$\pi$:  $\Bbb R^{n+1}$ $\x$ $\Bbb R$ $\to$ $\Bbb R^{n+1}\x\{0\}$ 
be the natural projection map.    
We suppose $\pi\vert_K$ is a self-transverse immersion.  
The {\it projection} $P$ of an $n$-knot $K$ is 
$\pi\vert_K(K)$ of $\Bbb R^{n+1}$.    
We give $P$ an orientation by using the orientation of $K$ naturally.  
The {\it singular point set} of the projection of an $n$-knot $K$ is 
the set  $\{x\in$ $\pi\vert_{K}$ $(K)$  $\vert$  
$\sharp$ $\{(\pi\vert_{K})$ $^{-1}(x)\}$$\geqq2.\}$.     
Let $\mu(P)$ denote the number of the connected components of 
the singular point set of the projection $P$.

Let $K$ be an $n$-knot with a projection $P$. 
Then the number $\mu(P)$ measures the complexity of $K$ as follows. 

Let $n$=1. 
If  $\mu(P)$ $\leqq2$, then $K$ is unknotted.  
(It is proved by chcking all posibble projections concretely. )

Let $n$=2. Suppose the singular point set of $P$ consists of double points.  
If  $\mu(P)$ $\leqq2$, then $K$ is unknotted.

Let $n$ be any natural number. 
There is an $n$-knot $K$ with a projection $P$  with the following properties.  
(1)$\mu(K)$=3 
(2)$K$ is truly knotted.      
(3)The singular point set of $P$ consists of double points.  
 Proof. Let $K_1$ be the trefoil knot. 
 Let $K_n$ be the (0-twist) spun knot of $K_{n-1}$ ($n\geqq2$). 
( See [Z] for  twist spun knots. )

It is natural to consider the following problem. 

\definition { Problem A }
Let  $K$ be  an $n$-knot with a Projection $P$
(thus the underlying manifold $K$ is an $n$-sphere ).   
Suppose the singular point set of $P$ consists of double points.  
Suppose that $\mu(P)\leqq2$.    
Then, is  $K$ unknotted?     
\enddefinition

Of course, if $n=1$ or $2$, as mentioned above, then 
the answer is affirmative.  
But for general $n\geq5$, we prove 
the answer to Problem A is negative in \S3. 

 We prove:
\proclaim 
{Theorem 1}
Let $n\geqq5$.
There is an $n$-knot $K$ with a projection $P$ with the following properties. 

(1)  $K$ is diffeomorphic to the standard $n$-sphere. 

(2)  The singular point set of $P$ consists of double points.  
   
(3)  $\mu(P)=2$.        

(4)  $K$ is truly knotted. 
\endproclaim 

In \S4,  furthermore, we prove  
the projection of the $n$-knot constructed in the proof of Theorem 1 
has the following property. 

\proclaim 
{Theorem 2}
Let $n\geqq5$.
There is an $n$-knot with a projection $P$ such that 
$P$ is not the projection of any knot which is unknotted. 
\endproclaim

\f{\bf Note.} 
(1)
It is well-known that the projection of any 1-dimensional knot 
is the projection of a 1-knot which is unknotted.   
The fact is used in definitions of the Jones polynomial 
and the Conway-Alexander polynomial.   
See [Kf1] and [Kf2].

(2)
The author proved the $n\geqq3$ case of Theorem 2 is true
in [O].  
But $\mu(P)$ of the examples are greater than two.

In the case of codimension two submanifolds of $\Bbb R^{n+2}$
which are diffeomorphic to a connected closed manifold and which are not spheres,  
we have the following Problem B corresponding to Problem A.

Let $M$ be a connected closed $n$-manifold.
Let  $K$ be a submanifold of $\Bbb R^{n+2}$ which is diffeomorphic to $M$.   
Suppose $\pi\vert K$ is transverse. 
Put $P=\pi(K)$. 
The number $\mu(P)$ is defined similarly. 

Submanifolds $K_1$ and $K_2$ in $\Bbb R^{n+2}$ are said to be {\it equivalent} 
if there is a diffeomorphism $f:\Bbb R^{n+2}\rightarrow\Bbb R^{n+2}$ 
such that $f(K_1)=K_2$ and that $f\vert K_1$ is an orientation preserving diffeomorphism 
if $M$ is oriented.

\definition{Problem B}
Let $M$, $K$, $P$ and $\mu(\quad)$ be as above.  
Suppose the singular point set of $P$ consists of double points.  
Suppose that $\mu(P)$ $\leqq2$ (resp. $\leqq1$).   
Then, is an equivalence class of submanifolds determined uniquely? 
In particular, is it determined uniquely when $M$ is embedded in $\Bbb R^{n+1}$?  
\enddefinition  

[Sh] shows that, when $M$ $\cong$ $T^2$, 
then the equivalence class of submanifolds is determined provided 
$\mu(\quad)$ $\leqq2$.  
On the other hand, for high dimensional case we have the following.

In \S2 we prove:        

\proclaim{Theorem 3 }    
Let $n\geqq5$. 
There is a closed connected oriented  $n$-dimensional manifold $M$ as follows. 
There are submanifolds $K_i$ with a projection $P_i$ ($i=0,1$) 
which are diffeomorphic to $M$ with the following properties. 

(1) 
$\mu(P_0)$=0. 

(2) 
$\mu(P_1)$=1.

(3)  The singular point set of $P_i$ consists of double points.  

(4)
$K_0$ is equivalent to neither $K_1$, $-K_1$, $K_1^*$ nor $-K_1^*$. 

(5)
$M$ is embedded in $\Bbb R^{n+1}$.
\endproclaim

The construction of the manifold $M$ in Theorem 3 will be used 
in the proofs of Theorem 1 and 2.

\head 2. The proof of Theorem 3 \endhead

We first prove the case of $n=5$. 

We define submanifolds 
$K_0$ and $K_1$ $\subset$ $\Bbb R^7$=$\Bbb R^6\x\Bbb R^1$ 
which are diffeomorphic to  $S^3\x S^2$.
Of course $S^3\x S^2$ is embedded in $\Bbb R^6$.

We  define $K_0$ $\subset$ $\Bbb R^7$=$\Bbb R^6\x\Bbb R^1$.  
Let $A_0$ be a trivially  embedded $3$-sphere in $\Bbb R^6$$\x\{0\}$. 
Take the tubular neighborhood $N_0$ of $A_0$ in $\Bbb R^6$$\x\{0\}$.
Then $\partial N_0$ is diffeomorphic to $S^3\x S^2$. 
Define  $K_0$ to be  $\partial N_0$.  
The projection $P_0$ of $K_0$ coincides with $K_0$. 
Obviously $\mu(P_0)$=0.

We define $K_1$ $\subset$ $\Bbb R^7$ =$\Bbb R^6\x\Bbb R^1$.  
Take a self-transverse immersion 
$g:S^3\i$ $\Bbb R^6\x\{0\}$ 
such that the singular point set is one point $p$.  
Then $\sharp\{g^{-1}(p)\}=2$. 
We suppose that there is a subset $V$ of $\Bbb R^6\x\Bbb R^1$  with the following properties. 

\roster
\item
$V=\{(x_1,x_2,x_3,y_1,y_2,y_3,z)\vert$ 
$x_1^2+x_2^2+x_3^2$ $<1$,  
$y_1^2+y_2^2+y_3^2$ $<1$, 
$z\in\Bbb R.\}$. 

\item
$V\cap $$g(S^{3})$ is a union of two open $3$-discs 
$D^3_x$ and $D^3_y$.

\item
$D^3_x=\{(x_1,x_2,x_3,y_1,y_2,y_3,z)\vert$ 
$x_1^2+x_2^2+x_3^2$ $<1$,  
$y_1$=$y_2$=$y_3$=0, $z=0$ $ \}$. 

\item
$D^3_y=\{(x_1,x_2,x_3,y_1,y_2,y_3,z)\vert$ 
$x_1$=$x_2$=$x_3$=0, 
$y_1^2+y_2^2+y_3^2$ $<1$, 
$z=0$ 
$ \}$. 
\endroster

Take the normal bundle $\nu$ of $g(S^{3})$ in $\Bbb R^{6}\x\{0\}$. 
Let $E$ be a manifold which is the total space of $\nu$. 
Thus we obtain an immersion  
$\widetilde{g}:$ $E\i$ $\Bbb R^{6}\x\{0\}$. 
Since $\pi_2 SO(3)$=0, $\nu$ is the trivial bundle  
and 
$\partial E$ is diffeomorphic to   $S^3\x S^2$.

We can take $\widetilde{g}$ to satisfy the following conditions. 

\roster
\item
 $\widetilde{g}\vert_{V^C}$ is an embedding,  
where  $V^C$ is  $\widetilde{g}^{-1}$ $(\widetilde{g}(E)-\{\widetilde{g}(E)\cap V\})$. 
 \item   
$\widetilde{g}(E)\cap V$  
\newline  
= $\{(x_1,x_2,x_3,y_1,y_2,y_3,z)\vert$ 
$x_1^2+x_2^2+x_3^2$ $<1$,  
$y_1^2+y_2^2+y_3^2$ $\leqq$ $\frac{1}{4}$, $z=0\}$  

$\cup$ 
$\{(x_1, x_2, x_3,y_1, y_2, y_3, z)\vert$ $x_1^2+x_2^2+x_3^2$ $\leqq$ $\frac{1}{4}$, 
$y_1^2+y_2^2+y_3^2$ $<1$,  $z=0\}$. 

\endroster



\vskip1cm
 Figure 1.

You can obtain this figure 
by clicking `PostScript' in the right side of 
the cite of the abstract of this paper in arXiv 
(https://arxiv.org/abs/the number of this paper). 

You can also obtain it from the author's website,  
which can be found by typing his name in search engine. 

\vskip1cm

Let $f:E\e$ $\Bbb R^6\x\Bbb R^1$ 
be an embedding with the following properties.

\roster
\item
 $\widetilde{g}\vert_{V^C}$=   $f\vert_{V^C}$   

\item   
$f(E)\cap V$  
\newline  
=         
$\{(x_1, x_2, x_3,y_1, y_2, y_3, z)\vert$ 
$x_1^2+x_2^2+x_3^2$ $<1$,  
$y_1^2+y_2^2+y_3^2$ $\leqq$ $\frac{1}{4}$, 

$z=$ $1-(x_1^2+x_2^2+x_3^2)$
$\}$  

$\cup$ 
$\{(x_1, x_2, x_3,y_1, y_2, y_3, z)\vert$ 
$x_1^2+x_2^2+x_3^2$ $\leqq$ $\frac{1}{4}$, 
$y_1^2+y_2^2+y_3^2$ $<1$,  
$z=0$
$\}$. 
\endroster
We can make the corner smooth.  


\vskip1cm
 Figure 2. 

You can obtain this figure 
by clicking `PostScript' in the right side of 
the cite of the abstract of this paper in arXiv 
(https://arxiv.org/abs/the number of this paper). 

You can also obtain it from the author's website,  
which can be found by typing his name in search engine. 

\vskip1cm

The submanifold $f(\partial E)$ is called $K_1$.  
Then the projection $P_1$ of $K_1$ is  $\widetilde{g}(\partial E)$.

Then we have: 
The singular point set of the projection $P_1$  is 
 \nl
$\{(x_1, x_2, x_3, y_1, y_2, y_3, z)\vert$ 
$x_1^2+x_2^2+x_3^2$=$\frac{1}{4}$, 
$y_1^2+y_2^2+y_3^2$=$\frac{1}{4}$, $z=0\}$. 
It consists of double points.  
It is diffeomorphic to 
$S^{2}\x S^{2}$.   
It is  connected. 
$\mu(P_1)=1$.

We prove:  $K_0$ is not equivalent to any of $K_1$, $-K_1$, $K_1^*$ nor $-K_1^*$.


Proof.
Let $K$ be a codimension two submanifold of $S^{n+2}$. 
Let $X_K$ denote the infinite cyclic covering space of the complement 
associated with the natural homomorphism map 
$\pi_1(S^{n+2}-K)\rightarrow H_1(S^{n+2}-K;\Bbb Z)\cong\Bbb Z$. 
We consider $H_*(X_K;\Bbb Z)$ as a module over $\Lambda=\Bbb Z[t,t^{-1}]$.  
See [M], [L1] etc. for properties of such spaces and those of such modules.  

We can regard that $K_i$ is in $S^7$($i=0,1$) naturally.                              
We consider $H_3(X_{K_i};\Bbb Z)$. 
By the construction of $K_i$, we have:

(1) $H_3(X_{K_0};\Bbb Z)\cong0$. 

(2) $H_3(X_{K_1};\Bbb Z)\cong$
 $H_3(X_{-K_1};\Bbb Z)\cong$
 $H_3(X_{K_1^*};\Bbb Z)\cong$
 $H_3(X_{-K_1^*};\Bbb Z)\cong\Lambda/(t-1)\cdot\Lambda$.

Therefore $K_0$ is equivalent to neither $K_1$, $-K_1$, $K_1^*$ nor $-K_1^*$.

We next prove the case of $n>5$. 

We define an $n$-dimensional submanifold  $K_i^{(n)}\subset\Bbb R^{n+2}$ 
as follows.    ( $n\geqq5$, $i=0,1$.) 

Let $K_i^{(5)}$ be $K_i$. 

Put $\Bbb R^{n+2}=\{x \vert x\in\Bbb R\}\x$ $\Bbb R^n\x$ 
$\{t\vert t\in\Bbb  R\}$. 
Suppose the projection map is 
$\Bbb R^{n+2}\rightarrow$  $\{x\vert x\in\Bbb R\}\x$ $\Bbb R^n$$\x$ $\{t\vert t=0\}$.

We assume 
$K_i^{(n)}\subset$$\{x\vert x\geqq0\}\x$ $\Bbb R^n$$\x$ $\{t\vert t\in\Bbb R\}$  
and 
$K_i^{(n)}\cap$$\{x\vert x=0\}\x$ $\Bbb R^n$$\x$ $\{t\vert t\in\Bbb R\}$ 
is an $n$-disc. 

We define $K_i^{(n+1)}\subset\Bbb R^{(n+3)}$ as follows.  
We consider 
$\Bbb R^{(n+3)}=\{(x,y)\vert x,y\in\Bbb R\}\x$ $\Bbb R^n$$\x$ $\{t\vert t\in\Bbb R\}$. 
We regard $\Bbb R^{(n+3)}$ as the result of rotating  
$\{x\vert x\geqq0\}\x$ $\Bbb R^n$$\x$ $\{t\vert t\in\Bbb R\}$ 
around 
$\{x\vert x=0\}\x$ $\Bbb R^n$$\x$ $\{t\vert t\in\Bbb R\}$. 
When rotating it, rotate $K_i^{(n)}$ as well. 
The result is called $K_i^{(n+1)}$.

By the construction of $K_i^{(n)}$, we have: 

(1) For the projection of $P_i^{(n)}$ of $K_i^{(n)}$,   $\mu(P_i^{(n)})$=1. 
    The singular point set of $P_i$ consists of double points.  

(2) $H_3(X_{K^{(n)}_0};\Bbb Z)\cong0$. 

(3) $H_3(X_{K^{(n)}_1};\Bbb Z)\cong$
 $H_3(X_{-K^{(n)}_1};\Bbb Z)\cong$
 $H_3(X_{K^{(n)*}_1};\Bbb Z)\cong$
 $H_3(X_{-K^{(n)*}_1};\Bbb Z)\cong\Lambda/(t-1)\cdot\Lambda$.

The computation for $K_1$ follows because the 0-section of $E$ is a generator for 
$H_3(E)$ (Compare [Kf1], p.43, 190, 229).

Therefore $K^{(n)}_0$ is equivalent to neither 
$K^{(n)}_1$, $-K^{(n)}_1$, $K^{(n)*}_1$ nor $-K^{(n)*}_1$.

\head 3. 
The proof of Theorem 1
\endhead

We first prove the case of $n=5$. 

We use $f(E)$ in \S2. 

We suppose that $f(E)-V\subset$ $\Bbb R^6\x\{0\}$.
Take a 6-ball $B^6$ $\subset$ $\Bbb R^6\x\{0\}$ 
$\subset$ $\Bbb R^6\x\Bbb R$.  
In $B^6\x\Bbb R$, take a submanifold $A_1$ 
which is a parallel displacement of the submanifold $f(E)$.  
In $(\Bbb R^6-B^6)$$\x\Bbb R$, 
take a submanifold $A_2$ 
which is a parallel displacement of the submanifold $f(E)$ 
with the opposite orientation.   

Recall $E=S^3\x D^3$.  
We can put $A_i=S^3_i\x D^3$=
$(D^3_{iS}\cup D^3_{iN})$ $\x D^3$=
$(D^3_{iS}\x D^3)\cup(D^3_{iN}\x D^3)$ (i=1,2). 
Suppose 
$(D^3_{iS}\x D^3)$ is embedded in $\Bbb R^6\x\{0\}$.

Take submanifolds $S^2_1$ and $S^2_2$ diffeomorphic to the 2-sphere 
in $\partial B^6$    so that the linking number is one.


\vskip1cm
 Figure 3.

You can obtain this figure 
by clicking `PostScript' in the right side of 
the cite of the abstract of this paper in arXiv 
(https://arxiv.org/abs/the number of this paper). 

You can also obtain it from the author's website,  
which can be found by typing his name in search engine. 

\vskip1cm

There are orientation preserving diffeomorphism maps $h$ such that 
$h(\Bbb R^6\x\{t\})$=$\Bbb R^6\x\{t\}$.

By using such a diffeomorphism map, we move $A_1$ so that: 

(1) (Int$D^3_{1S})\x D^3$ $\subset(\Bbb R^6\x\{0\}-B^6)$
     
(2)(Int$D^3_{1N})\x D^3$ $\subset B^6\x\Bbb R$.   
 The singular point set of the projection of $A_1$ is in $B^6$. 

(3)$A_1\cap\partial B^6$=$(\partial D^3_{1S})\x D^3$
=$(\partial D^3_{1N})\x D^3$  and  
$\partial D^3_{1s}=\partial D^3_{1N}=S^2_1$.

By using such a diffeomorphism map, we move $A_2$ so that: 

(1) (Int$D^3_{2S})\x D^3$ $\subset B^6$

(2)(Int$D^3_{2N})\x D^3$ $\subset(\Bbb R^6-B^6)\x\Bbb R$.  
The singular point set of the projection of $A_2$ is in $\Bbb R^6-B^6$. 

(3)$A_2\cap\partial B^6$=$(\partial D^3_{2S})\x D^3$
=$(\partial D^3_{2N})\x D^3$ and  
$\partial D^3_{2s}=\partial D^3_{2N}=S^2_2$. 

We define $K$ to be 
\nl$\overline{\partial B^6-\partial (D^3_{1N}\x D^3)-\partial (D^3_{2N}\x D^3)}\cup$
$\overline{[\partial (D^3_{1N}\x D^3)\cup\partial (D^3_{2N}\x D^3)]-\partial B^6}$.


\vskip1cm
 Figure 4.

You can obtain this figure 
by clicking `PostScript' in the right side of 
the cite of the abstract of this paper in arXiv 
(https://arxiv.org/abs/the number of this paper). 

You can also obtain it from the author's website,  
which can be found by typing his name in search engine. 

\vskip1cm

By the construction, we have:

(1)  $K$ is  diffeomorphic to the 5-sphere.  $K$ is a 5-knot. 

(2) For the projection $P$ of $K$, $\mu(P)$=2.   
 The singular point set of $P$ consists of double points.

(3) A Seifert matrix of $K$ is 
$\pmatrix
1&1\\
0&-1\\
\endpmatrix.$
(See [L1] and [L2] for Seifert matrices.)
Hence $H_3(X_{K};\Bbb Z)\cong\Lambda/(t^2-3t+1)\cdot\Lambda$.  
Therefore $K$ is truly knotted. 

This completes the proof in the case of $n=5$. 

We next prove the case of $n>5$. 

Let $K^{(5)}$ be $K$. 
Let $K^{(n+1)}$ be the spun knot of $K^{(n)}$($n\geqq5$). 
(See [Z] for spun knots.)

We take the axis as in the proof of the $n>5$ case in \S 2. 
Then the projection $P^{(n+1)}$ of $K^{(n+1)}$is the result of 
rotating $P^{(n)}$ around the axis. 
Hence  $\mu(P)$=2.
The singular point set of $P$ consists of double points.

By the construction,  we have 
 $H_3(X_{K^{(n+1)}};\Bbb Z)\cong\Lambda/(t^2-3t+1)\cdot\Lambda$.  
Therefore $K^{(n+1)}$ is truly knotted. 

This completes the proof.

\head 4. 
The proof of Theorem 2
\endhead

We use $K$, $K^{(n)}$, and  $P^{(n)}$ in \S3. 

We first prove the case of $n=5$. 

Let $K'$ be a 5-knot. 
Suppose that the projection of $K'$ is the projection $P$ of $K$.  

Then a Seifert matrix of $K'$ is one of the following. 
\nl$\pmatrix
1&1\\
0&1\\
\endpmatrix$, 
$\pmatrix
1&1\\
0&-1\\
\endpmatrix$,  
$\pmatrix
1&-1\\
0&1\\
\endpmatrix$,   
$\pmatrix
1&-1\\
0&-1\\
\endpmatrix$, 
$\pmatrix
-1&1\\
0&1\\
\endpmatrix$, 
$\pmatrix
-1&1\\
0&-1\\
\endpmatrix$,  
$\pmatrix
-1&-1\\
0&1\\
\endpmatrix$,   
$\pmatrix
-1&-1\\
0&-1\\
\endpmatrix$. 
 
Hence  $H_3(X_{K'};\Bbb Z)$ is not trivial.  
Therefore $K'$ is truly knotted. 

We next prove the case of $n>5$. 

Let $K^{(n)'}$ be an $n$-knot($n>5$).   
Suppose the projection of  $K^{(n)'}$ is the projection $P^{(n)}$ of $K^{(n)}$.  
Then $K^{(n)'}$ is a spun knot of an $(n-1)$-knot 
whose projection is $P^{(n-1)}$. 
Hence  $H_3(X_{K^{(n)'}};\Bbb Z)$ is not trivial.  
Therefore $K^{(n)'}$ is truly knotted.

\np

\enddocument